\documentclass[12pt]{amsart}
\usepackage{amsfonts}
\newtheorem{thm}{Theorem}[section]
\newtheorem{theorem}[thm]{Theorem}

\newtheorem{lemma}[thm]{Lemma}

\newtheorem{example}[thm]{Example}
\newtheorem{remark}[thm]{Remark}

\theoremstyle{definition}
\newtheorem{definition}[thm]{Definition}

\numberwithin{equation}{section}

\newcommand{\dom}{\operatorname{dom}}
\newcommand{\slim}{\operatornamewithlimits{s-lim}}
\newcommand{\ov}{\overline}
\newcommand{\ti}{\widetilde}
\renewcommand{\Re}{\mathrm{Re}\,}

\newcommand{\bC}{{\mathbb C}}
\newcommand{\bN}{{\mathbb N}}
\newcommand{\bR}{{\mathbb R}}
\newcommand{\bZ}{{\mathbb Z}}

\newcommand{\fG}{{\mathfrak G}}
\newcommand{\fL}{{\mathfrak L}}
\newcommand{\fS}{{\mathfrak S}}

\newcommand{\al}{\alpha}
\newcommand{\be}{\beta}
\newcommand{\ga}{\gamma}
\newcommand{\La}{\Lambda}
\newcommand{\la}{\lambda}
\newcommand{\si}{\sigma}
\newcommand\om{\omega}

\begin{document}

\title[Half-inverse spectral problems]%
{Half-inverse spectral problems for Sturm--Liouville operators
    with singular potentials${}^{\dag}$}
\author[R.~O.~Hryniv, Ya.~V.~Mykytyuk]{Rostyslav O.~Hryniv and Yaroslav V.~Mykytyuk}%
\address{Institute for Applied Problems of Mechanics and Mathematics,
3b~Naukova st., 79601 Lviv, Ukraine and Lviv National University, 1
Universytetska st.,
79602 Lviv, Ukraine}%
\email{rhryniv@iapmm.lviv.ua, yamykytyuk@yahoo.com}%

\address{\emph{Current address of R.H.:}
Institut f\"ur Angewandte Mathematik, Abteilung f\"ur
Wahr\-schein\-lich\-keitstheorie und Mathematische Statistik,
Wegelerstr.~6, D-53115 Bonn, Germany}
\email{rhryniv@wiener.iam.uni-bonn.de}

\thanks{${}^{\dag}$The work was partially supported by Ukrainian Foundation
for Basic Research, grant No.~01.07/00172.
R.~H. acknowledges the support of the Alexander von Humboldt foundation}%


\date{\today}%

\begin{abstract}
Half-inverse spectral problem for a Sturm--Liouville operator consists in
reconstruction of this operator by its spectrum and half of the potential.
We give the necessary and sufficient conditions for solvability of the
half-inverse spectral problem for the class of Sturm--Liouville operators
with singular potentials from the space $W^{-1}_2(0,1)$ and provide the
reconstruction algorithm.
\end{abstract}
\maketitle


\section{Introduction}\label{sec:intro}

Assume that a function $q$ is integrable and real-valued on~$(0,1)$ and
that $h_0$ and $h_1$ are some elements of $\ov{\bR}:=\bR\cup\{\infty\}$.
Denote by $T=T(q,h_0,h_1)$ a Sturm--Liouville operator in $L_2(0,1)$ that
is given by the differential expression
\begin{equation}\label{eq:1.de}
    \ell:= -\frac{d^2}{dx^2} + q
\end{equation}
and the boundary conditions
\begin{equation}\label{eq:1.bc}
    y'(0) - h_0 y(0) = 0, \qquad y'(1) - h_1 y(1) = 0
\end{equation}
(where, as usual, $h_0=\infty$ or $h_1=\infty$ is a shorthand notation for
the Dirichlet boundary condition at the point $x=0$ or $x=1$ respectively).
Suppose that the spectrum $\Sigma(q,h_0,h_1)$ of~$T$, the number~$h_0$
and/or~$h_1$, and the potential~$q$ over a half of the interval---e.g.,
over~$(0,\tfrac12)$---are known; can one recover the operator~$T$ based on
this information? Problems of such kind are known in the literature as
half-inverse spectral problems, or inverse problems with mixed spectral
data.

The main aim of the present paper is to study the half-inverse spectral
problem for Sturm--Liouville operators with real-valued singular potentials
from the space $W^{-1}_2(0,1)$ (see Section~\ref{sec:form} for precise
definitions). Namely, we shall find necessary and sufficient conditions on
mixed data in order that the half-inverse spectral problem be soluble in
the considered class of potentials. We also specify this result to the case
of regular potentials from $L_2(0,1)$ and establish then a local existence
theorem.

The first result on the half-inverse spectral problem is due to Hochstadt
and Lieberman~\cite{HL}, who proved that if $\Sigma(\,\ti
{\vphantom{h}q},\ti h_0,\ti h_1) =\Sigma(q,h_0,h_1)$, $\ti h_0 = h_0$, $\ti
h_1=h_1$, and $\ti q = q$ on $(0,\tfrac12)$, then $\ti q = q$ on $(0,1)$.
Later, Hald~\cite{Ha1} proved that the statement remains true even if the
boundary conditions at the point $x=1$ are not assumed equal, while
del~Rio~\cite{dR} constructed counterexamples demonstrating that uniqueness
might fail if $\ti h_0 \ne h_0$. Hald~\cite{Ha} generalized the theorem by
Hochstadt and Lieberman to the setting motivated by the inverse problem for
the torsional modes of the Earth, where the domain of $T$ forces a
discontinuity in an interior point. It was shown in~\cite{Ha} that not only
the potential in $(\tfrac12,1)$, but also the position of the discontinuity
and the jump magnitude are uniquely determined by the spectrum
$\Sigma(q,h_0,h_1)$, the number $h_0$, and the potential $q$ on
$(0,\tfrac12)$. Willis~\cite{Wi} established similar uniqueness result for
an inverse problem with two interior discontinuities, and
Kobayashi~\cite{Ko} proved that a symmetric potential $q\in L_1(0,1)$ with
two symmetrically placed jumps is uniquely determined from $h_0$ and a
single spectrum $\Sigma(q,h_0,-h_0)$. We also mention the paper~\cite{bK},
where the uniqueness for the half-inverse problem was established for the
operator generated by the differential expression
 \[
    \frac1{a(x)} \frac{d}{dx}\Bigl( a(x) \frac{d}{dx} \Bigr) + q(x)
 \]
with the impedance $a$ vanishing at $x=0$ and positive otherwise.

Afterwards, these uniqueness results have been further generalized to
different settings. In~\cite{dRGS}, del Rio, Gesztesy and Simon proved a
number of results when $q$ on a part of $(0,1)$ and certain parts of the
spectra $\Sigma(q,h_0,h_1)$ for several different values of $h_1$
completely determine $q$ on $(0,1)$---e.g., so do the spectrum
$\Sigma(q,h_0,h_1)$, half the spectrum $\Sigma(q,h_0,h'_1)$ for $h_1'\ne
h_1$, and $q$ on $(0,\tfrac14)$, or two-thirds of the spectra
$\Sigma(q,h_0,h_1)$ for three different values of $h_1$. In~\cite{GS},
Gesztesy and Simon proved that it suffices to know all the eigenvalues of
$T$ except $k+1$ provided the potential $q$ is $C^{2k}$-smooth at the
midpoint $x=\tfrac12$. It was also proved in~\cite{GS} that, under suitable
growth conditions on the potential~$q$ of a Schr\"odinger operator~$T$ on
the whole line, the spectrum of $T$ and the potential on the half-line
$(0,\infty)$ determine uniquely $q$ on the other half-line $(-\infty,0)$.
This last result was recently improved by Khodakovsky~\cite{Kh}.

 In the recent work~\cite{Sakh}, Sakhnovich studied existence of
solution to the half-inverse spectral problem; namely, he presented
sufficient conditions on a function $q_0$ on $(0,\tfrac12)$ and a sequence
$(\la_n^2)$ of pairwise distinct real numbers tending to $+\infty$ in order
that there existed a Sturm--Liouville operator $T$ on~$(0,1)$ with given
spectrum $\{\la_n^2\}$ and potential $q$ coinciding with $q_0$ on
$(0,\tfrac12)$. These sufficient conditions are of local nature as they
require that $q_0$ belong to~$\mathrm{C}^1[0,\tfrac12]$ and that some
quantities constructed through the mixed spectral data $\{(\la_n^2),q_0\}$
be small enough; basically, $q_0$ is required to have small enough norm in
$\mathrm{C}^1[0,\tfrac12]$ and the sequence $(\la_n^2)$ to be close enough
to the ``unperturbed'' one, $(\pi^2n^2)$. Under these smallness assumptions
the author also suggested a constructive algorithm based on iterative
solution of a Gelfand--Levitan--Marchenko type integral equation.

We also mention the book by P\"oschel and Trubowitz~\cite{PT}, where among
other inverse results the authors prove that the odd part of a potential
$q\in L_2(0,1)$ and the Dirichlet spectrum $\Sigma(q,\infty,\infty)$
determine uniquely the whole potential. Coleman and McLaughlin~\cite{CM}
generalized the approach by P\"oschel and Trubowitz and derived an
analogous result for Sturm--Liouville operators in impedance form
\[
    Su = \frac1{p(x)}\frac{d}{dx}\Bigl( p(x) \frac{d}{dx} \Bigr)
\]
with $p\in W^1_2(0,1)$. Here $p$ is a positive \emph{impedance}, and the
operator $S$, subject to suitable boundary conditions, is selfadjoint in
the weighted space $L_2((0,1);p\,dx)$.

Sturm--Liouville operators in impedance form were earlier treated by
Andersson~\cite{An1,An2}.  In \cite{An1} the cases $p \in W^1_r(0,1)$,
$r\in[1,\infty)$, and $p$ of bounded variation were considered and various
direct and inverse spectral problems addressed; in particular, uniqueness
and local existence results for the half-inverse spectral problem were
derived. In~\cite{An2}, among other direct and inverse spectral problems, a
global existence result for the half-inverse spectral problem for the
impedance Sturm--Liouville operator $S$ with $p\in W^1_2(0,1)$ was claimed.
Namely, Theorem~5.2 of~\cite{An2} states that for any positive function
$p_0\in W^1_2(0,\tfrac12)$ and any sequence $(\la^2_n)$ of pairwise
distinct real numbers tending to~$+\infty$ and obeying the necessary
asymptotics there exists an impedance Sturm--Liouville operator $S$ with
impedance $p\in W^1_2(0,1)$ extending $p_0$, whose spectrum coincides with
the set $\{\la_n^2\}$. Unfortunately, as the example of
Section~\ref{sec:main} demonstrates, Theorem~5.2 of~\cite{An2} is
erroneous.

We observe that for $p\in W^1_2(0,1)$ the impedance
Sturm--Liouville operator~$S$ is unitarily equivalent to a
standard Sturm--Liouville operator with singular potential
$q:=(\sqrt{p})''/\sqrt{p}\in W^{-1}_2(0,1)$. The theory of
Sturm--Liouville operators with singular potentials
from~$W^{-1}_2(0,1)$ have been thoroughly developed in the recent
works by Shkalikov and Savchuk (see, e.g., \cite{SS,SS1}), and
some settings of inverse spectral problems for such operators have
been completely solved in~\cite{HMinv,HMinv2,HMinv3}.

The present work was highly motivated by the papers~\cite{An2,Sakh}: we
wanted, firstly, to correct the erroneous statement of Theorem~5.2
in~\cite{An2} concerning solvability of the half-inverse spectral problem
and, secondly, to formulate local existence result analogous to that
of~\cite{Sakh} directly in terms of the mixed data.

The paper is organized as follows. In the next section we introduce
necessary definitions and formulate the main results. In
Section~\ref{sec:transf} some facts about transformation operators are
presented. Proof of Theorem~\ref{thm:main} and a counterexample are given
in Section~\ref{sec:main}, and Theorems~\ref{thm:local} and
\ref{thm:regular} are proved in Section~\ref{sec:prf2.2}. Finally,
Appendix~\ref{sec:smooth} contains the proof of Theorem~\ref{thm:3.smooth}.

%
\section{Formulation of the main results}\label{sec:form}
%

Throughout the paper we shall denote by $(\cdot,\cdot)$ and $\|\cdot\|$
respectively the scalar product and the norm in the Hilbert space
$L_2(0,1)$. The symbol $W_p^s(a,b)$, $p\ge1$, $s\in\bR$, shall stand for
the Sobolev function space over $(a,b)$, and we write $\Re L_2(a,b)$ and
$\Re W^1_2(a,b)$ for the real Hilbert spaces of real-valued functions from
$L_2(a,b)$ and $W^1_2(a,b)$ respectively.

Suppose that $q$ is a real-valued distribution from $W^{-1}_2(0,1)$ and
that $\si\in L_2(0,1)$ is any of its (real-valued) distributional
primitive. Then the differential expression $\ell$ of~\eqref{eq:1.de} can
be written as
\[
    \ell_\si(y) = - (y'-\si y)' - \si y'=: - \bigl(y^{[1]}\bigr)' - \si y'.
\]
In what follows, the symbol $y^{[1]}$ will stand for the
\emph{quasi-derivative} $y'- \si y$ of a function $y\in W^1_1(0,1)$. A
natural $L_2$-domain of $\ell_\si$ is
\[
    \dom \ell_\si = \{y\in W^1_1(0,1) \mid y^{[1]} \in W^1_1(0,1),
                l_\si(y)\in L_2(0,1)\}.
\]
For any $h_0,h_1\in \bR\cup \{\infty\}$ we denote by $T=T(\si,h_0,h_1)$ an
operator in $L_2(0,1)$ that acts as $Ty = \ell_\si(y)$ on the domain
\[
    \dom T = \{y\in\dom \ell_\si\mid y^{[1]}(0) = h_0 y(0),
        y^{[1]}(1) = h_1 y(1)\}.
\]

This method of regularization by quasi-derivatives was suggested by
Shkalikov and Savchuk in~\cite{SS} (see also \cite{SS1}). Since
$\ell_\si(y) = -y'' + \si' y$ in the sense of distributions, so defined
operator~$T$ can be regarded as a Sturm--Liouville operator with the
singular potential $q=\si'\in W^{-1}_2(0,1)$. It is known~\cite{SS} that
$T$ is a selfadjoint bounded below operator with discrete spectrum.

Note that $T(\si+h,h_0,h_1)=T(\si,h_0+h,h_1+h)$ for any $h\in\bR$;
henceforth, by adding a suitable constant $h$ to the primitive $\si$, we
can always achieve one of the following situations: (1)~$h_0=0$,
$h_1\in\bR$; (2)~$h_0=0$, $h_1=\infty$, (3)~$h_0=\infty$, $h_1=0$,
(4)~$h_0=h_1=\infty$. All four cases are treated in a similar way, and, to
be definite, we shall concentrate on the first case here.

Therefore we shall assume that the operator $T$ has the form
$T(\si,0,h)=:T_{\si,h}$ for some $\si\in \Re L_2(0,1)$ and some $h\in\bR$.
As was mentioned above, $T$ so defined is selfadjoint, bounded below, and
has simple discrete spectrum accumulating at $+\infty$. By adding a
suitable constant~$C$ to the potential~$q$ (or a suitable linear function
$Cx$ to the primitive~$\si$) if necessary, we can make all the eigenvalues
$\la^2_n$, $n\in\bZ_+$, of $T$ positive and will tacitly assume this in
what follows. It is known~\cite{Sav,HMtr} that the eigenvalues~$\la_n^2$,
when arranged in increasing order, obey the following asymptotic formula:
\[
    \la_n = \pi n + \mu_n, \qquad n\in\bZ_+,
\]
where $(\mu_n)_{n\in\bZ_+}\in \ell_2$.

We denote by $\fL$ the set of all strictly increasing sequences
$\La=(\la_n)_{n\in\bZ_+}$, in which $\la_n$ are positive numbers such that
$\mu_n:=\la_n-\pi n$ form an $\ell_2$-sequence. Fix an arbitrary
$\La=(\la_n)\in\fL$ and denote by $\Pi_\La$ the set of all real-valued
functions $\psi\in L_2(0,1)$ of the form
\begin{equation}\label{eq:2.PiLa}
    \psi(x) = \sum_{n=0}^\infty[\al_n \cos(\la_nx) - \cos (\pi nx)]+\tfrac12,
\end{equation}
where $(\al_n)_{n\in\bZ_+}$ is a sequence of positive numbers such that the
sequence $(\al_n-1)_{n\in\bZ_+}$ belongs to $\ell_2$. Put
\begin{equation}\label{eq:2.psiLa}
    \psi_\La(x) := \sum_{n=0}^\infty[\cos(\la_nx) - \cos (\pi nx)];
\end{equation}
it can be shown (see Lemma~\ref{lem:4.psiLa}) that this series
converges in~$L_2(0,1)$. Since the system
$\{\cos(\la_nx)\}_{n=0}^\infty$ forms a Riesz basis of
$L_2(0,1)$~\cite{HV}, we conclude that for any $\psi\in L_2(0,1)$
there is an $\ell_2$-sequence $(\be_n)$ such that
\[
    \psi(x)-\psi_\La(x)-\tfrac12 = \sum_{n=0}^\infty \be_n\cos\la_nx.
\]
It follows that $\psi$ admits representation~\eqref{eq:2.PiLa} with
$\al_n:=\be_n+1$. Thus the only restriction imposed by $\Pi_\La$ is that
all the coefficients $\be_n$ in the above series representation of
$\psi-\psi_\La-\tfrac12$ should be greater than $-1$. The above arguments
also show that $\Pi_\La$ is an open and convex set in $\Re L_2(0,1)$.

Assume that $\si_0 \in \Re L_2(0,\tfrac12)$ and denote by
$y_0(\cdot,\la)=y_0(\cdot,\la,\si_0)$, $\la\in\bC$, the solution of the
equation
\[
    -(y'-\si_0y)' - \si_0 y' = \la^2 y
\]
on the interval $(0,\tfrac12)$ subject to the initial conditions
\[
    y(0) = 1,\quad (y'- \si_0 y)(0)=0.
\]
Let $l_{0,\si_0}(x,t)$ be the kernel of the transformation operator
$I+L_{0,\si_0}$ that maps $y_0(\cdot,\la)$ into the function $\cos\la x$
for all complex $\la$ (see Section~\ref{sec:transf} for details). In other
words, $l_{0,\si_0}$ is such that the following equality is satisfied for
all $\la\in\bC$, $x\in(0,\tfrac12)$:
\[
    \cos\la x = y_0(x,\la) + \int_0^x l_{0,\si_0}(x,t)y_0(t,\la)\,dt.
\]
We put
\begin{equation}\label{eq:2.phi0}
    \phi_0(2x)=\phi_0(2x,\si_0) := - \frac12\si_0(x) + \int_0^{x}
    l^{\,2}_{0,\si_0}(x,t)\,dt, \quad x\in(0,\tfrac12).
\end{equation}
It follows from the results of~\cite{HMtr} that $\phi_0\in L_2(0,1)$; see
also Remark~\ref{rem:3.coinc}.

Our main result is the following
\begin{theorem}\label{thm:main}
Assume that $\La=(\la_n)_{n\in\bZ_+}\in\fL$, $\si_0\in \Re
L_2(0,\tfrac12)$, and define $\phi_0:=\phi_0(\cdot,\si_0)$ as
in~\eqref{eq:2.phi0}.
\begin{itemize}
\item[(i)] The half-inverse spectral problem is soluble for the mixed
spectral data $\{\si_0,\La\}$ if and only if the function $\phi_0$ belongs
to $\Pi_\La$.
\item[(ii)] If $\phi_0\in \Pi_\La$, then the solution of the above
half-inverse spectral problem is unique---i.e., there exist a unique
$\si\in\Re L_2(0,1)$ and a unique $h\in\bR$ such that $\si$ is an extension
of $\si_0$ and the spectrum of $T_{\si,h}$ coincides with
$\La^2:=(\la_n^2)_{n\in\bZ_+}$.
\end{itemize}
\end{theorem}

As the function $\phi_0$ is determined via $\si_0$ and the set $\Pi_\La$
via $\La$, condition (i) of the theorem imposes connection between $\si_0$
and $\La$. An example of mixed data $\{\si_0,\La\}$, for which $\phi_0$ is
not in $\Pi_\La$ (and thus the half-inverse problem has no solution), is
given at the end of Section~\ref{sec:main}.

For the ``unperturbed'' situation with zero potential one has
$\si_0\equiv0$ on $(0,\tfrac12)$, $\la_n=\pi n$, $n\in\bZ_+$, and
$\phi_0\equiv\tfrac12$, so that $\phi_0\in\Pi_\La$. Since the set $\Pi_\La$
is open and depends continuously on $\La$, and the function $\phi_0$
of~\eqref{eq:2.phi0} depends continuously on~$\si_0$, it follows that
$\phi_0\in\Pi_\La$ if the function $\si_0\in\Re L_2(0,\tfrac12)$ and the
sequence $\La-(\pi n)\in\ell_2$ have small norms. This is a local existence
result analogous to those of the papers~\cite{An1,Sakh}. However, nice
bounds on the norms $\|\si_0\|_{L_2(0,1/2)}$ and $\|\La-(\pi n)\|_{\ell_2}$
are cumbersome and difficult to obtain. Instead, we shall establish such
bounds in the particular case where $\si_0\in W^1_2(0,\tfrac12)$, i.e.,
where $q_0\in L_2(0,\tfrac12)$, and estimate the norm of $q_0$ in
$L_2(0,\tfrac12)$.

\begin{theorem}\label{thm:local}
Assume that a real-valued function $q_0\in L_2(0,\tfrac12)$ and a sequence
$(\la_n)=:\La\in \fL$ are such that $\|q_0\|_{L_2(0,1/2)}\le\tfrac12$ and
$\|(\la_n-\pi n)\|_{\ell_2}\le\tfrac14$. Then $\phi_0\in \Pi_\La$;
therefore there exists a unique function $\si\in \Re L_2(0,1)$ extending
$\si_0:=\int_0^x q_0$ and a unique real~$h$ such that the numbers $\la_n^2$
are eigenvalues of the Sturm--Liouville operator~$T_{\si,h}$.
\end{theorem}

Theorems~\ref{thm:main} and \ref{thm:local} admit the following refinement
for the class of regular Sturm--Liouville operators with potentials
from~$L_2(0,1)$. We observe that if the potential $q$ belongs to $L_2(0,1)$
and $\si(x)=\int_0^x q(t)\,dt$, then the operator $T=T_{\si,h}$ is given by
\[
    Ty:= - y'' + qy, \qquad
    \dom T = \{y\in W^2_2(0,1) \mid y'(0)=0,\, y'(1) = \hat h y(1)\}
\]
with $\hat h:= h+ \si(1)$, and the eigenvalues $\la_n^2$ of $T$ obey the
asymptotics
\[
    \la_n = \pi n + \frac{c}n + \frac{\nu_n}{n}, \qquad n\in\bZ_+,
\]
for some $c\in\bR$ and an $\ell_2$-sequence $(\nu_n)$. We denote by $\fL_1$
a subset of $\fL$ formed by sequences $(\la_n)$ obeying this refined
asymptotics.

\begin{theorem}\label{thm:regular}
Assume that $\La\in\fL_1$ and $\si_0\in \Re W^1_2(0,\tfrac12)$. If
$\phi_0\in \Pi_\La$, then the extended function~$\si$ given by
Theorem~\ref{thm:main} belongs to $W^1_2(0,1)$.
\end{theorem}

In other words, this theorem states that if for a given function $q_0\in
\Re L_2(0,\tfrac12)$ and a given sequence $(\la_n)\in\fL_1$ the
half-inverse spectral problem has a solution within the class of
Sturm--Liouville operators with potentials from $W^{-1}_2(0,1)$, then the
recovered potential belongs in fact to $L_2(0,1)$.

\section{Transformation operators}\label{sec:transf}

In this section we shall formulate some results from the
papers~\cite{HMtr,HMinv} that will be used later on to establish our
principal results.

Suppose that $\si\in L_2(0,1)$ and denote by $\ti T_\si$ an operator in
$L_2(0,1)$ that acts according to
\[
    \ti T_\si y = \ell_\si(y) :=  - (y^{[1]})'-\si y'
\]
on the domain
\[
    \dom \ti T_\si = \{y\in \dom \ell_\si \mid y^{[1]}(0)=0\}.
\]
In other words, $\ti T_\si$ is an extension of the operator $T_{\si,h}$
discarding the boundary condition at the terminal point $x=1$.

One of the main results of the paper~\cite{HMtr} is that the operators $\ti
T_\si$ and $\ti T_0$ are similar, and the similarity is performed by a
\emph{transformation operator} of a special form.

\begin{theorem}\label{thm:3.tr}
Assume that $\si\in L_2(0,1)$; then there exists an integral
Hilbert--Schmidt operator $K_\si:\, L_2(0,1) \to L_2(0,1)$ of the form
\begin{equation}\label{eq:3.tr}
    (K_\si u)(x) = \int_0^x k_\si(x,t) u(t)\,dt
\end{equation}
such that $I+K_\si$ is a transformation operator for $\ti T_\si$ and $\ti
T_0$, i.e., such that
\begin{equation}\label{eq:3.sim}
    \ti T_\si (I + K_\si) = (I + K_\si) \ti T_0.
\end{equation}
The operator $K_\si$ with properties~\eqref{eq:3.tr}--\eqref{eq:3.sim} is
unique. If, moreover, the function $\si$ is real-valued, then the kernel
$k_\si$ is real-valued, too.
\end{theorem}

Put $L_\si:= (I+K_\si)^{-1} - I$; then $L_\si$ is an integral
Hilbert--Schmidt operator of Volterra type, i.e.,
\[
    (L_\si u)(x) = \int_0^x l_\si(x,t) u(t)\,dt,
\]
and $I+L_\si$ is the transformation operator for $\ti T_0$ and $\ti T_\si$.
The kernels $k_\si$ and $l_\si$ of the operators $K_\si$ and $L_\si$
possess the property that, for any fixed $x\in [0,1]$, the functions $k_\si
(x,\cdot)$ and $l_\si (x,\cdot)$ belong to $L_2(0,1)$ and the mappings
\begin{align*}
    [0,1]\ni x &\mapsto k_\si(x,\cdot) \in L_2(0,1),\\
    [0,1]\ni x &\mapsto l_\si(x,\cdot) \in L_2(0,1)
\end{align*}
are continuous.

The transformation operators naturally appear during factorization of some
Fredholm operators, which we shall now explain.

Denote by $\fS_2$ the ideal of all Hilbert--Schmidt operators in
$L_2(0,1)$. It is known that any operator in $\fS_2$ is an integral
operator with square integrable kernel on~$(0,1)\times(0,1)$. We denote by
$\fS^+_2$ ($\fS^-_2$) the subalgebra of $\fS_2$ consisting of all integral
operators in $\fS_2$ with upper-diagonal (respectively, lower-diagonal)
kernels.

\begin{definition}\label{def:3.fact}
    Assume that $Q\in \fS_2$. We say that the operator $I+Q$ is
    \emph{factorizable}, or that $I+Q$ \emph{admits factorization},
    if there exist operators $R^+\in\fS_2^+$ and $R^-\in\fS_2^-$
    such that
\[
    I+Q = (I+R^+)^{-1}(I+R^-)^{-1}.
\]
\end{definition}

Observe that an operator $I+Q$ can admit at most one factorization and thus
the operators $R^{\pm} = R^\pm(Q)$ are uniquely determined by $Q$.

To every function $\phi\in L_2(0,2)$ we shall put into correspondence an
integral operator $F_\phi\in \fS_2$ with kernel $f_{\phi}(x,t) := \phi(x+t)
+ \phi(|x-t|)$, i.e.,
\[
    (F_\phi u)(x) := \int_0^1 f_\phi(x,t) u(t)\,dt.
\]
Denote by $\Phi$ the set of those $\phi\in L_2(0,2)$, for which the
corresponding operator $F_\phi$ admits factorization. It follows from the
results of~\cite{Fact} that the set $\Phi$ is open and everywhere dense in
$L_2(0,2)$. We note also that if $F_\phi$ is a selfadjoint operator, then
$\phi\in \Phi$ if and only if $I+F_\phi$ is (uniformly) positive in
$L_2(0,1)$, see~\cite[Ch.~4]{GKvolt}.

Connection between the transformation operators $K_\si$ and operators
$F_\phi$ is described by the following statement, cf.~\cite{HMtr}.

\begin{theorem}\label{thm:3.conn}
\begin{itemize}
\item[(i)] Let $\si\in L_2(0,1)$ and define a function $\phi=\phi_\si$ via
\begin{equation}\label{eq:3.phi}
    \phi(2x) = - \frac12\si(x) + \int_0^x l^2_\si(x,t)\,dt, \qquad x\in(0,1).
\end{equation}
Then $\phi\in\Phi$ and
\begin{equation}\label{eq:3.fact}
    I+F_\phi = (I+K_\si)^{-1}(I+K_\si^\top)^{-1},
\end{equation}
where $I+K_\si$ is the transformation operator for $\ti T_\si$ and $\ti
T_0$, and $K_\si^\top$ is the operator associated to $K_\si$, i.e.,
\[
    (K_\si^\top u)(x):= \int_x^1 k_\si (t,x) u(t)\,dt.
\]
\item[(ii)] Conversely, if $\phi\in\Phi$ and $I+F_\phi = (I+K^+)^{-1}(I+K^-)^{-1}$
with $K^\pm\in\fS^\pm_2$, then $K=K_\si$ for some $\si\in L_2(0,1)$, $K^- =
K_\si^\top$, and~\eqref{eq:3.phi} holds.
\item[(iii)] The mapping $L_2(0,1)\ni \si \mapsto \phi \in \Phi$ given
by~\eqref{eq:3.phi} is homeomorphic.
\end{itemize}
\end{theorem}

The same statements hold certainly true if we consider the Sturm--Liouville
problem on the interval $(0,\tfrac12)$ instead of $(0,1)$; the function
$\phi$ will then be defined on $(0,1)$ instead of $(0,2)$.

It can be proved that for a smooth function~$\si$ the function
$\phi$ of~\eqref{eq:3.phi} is also smooth. We shall need the
following version of this statement (proved in
Appendix~\ref{sec:smooth}).

\begin{theorem}\label{thm:3.smooth}
The restriction of the mapping~$L_2(0,1)\ni \si \mapsto \phi \in \Phi$
given by~\eqref{eq:3.phi} to the set $W^1_2(0,1)$ is a bijection onto
$\Phi\cap W^1_2(0,2)$.
\end{theorem}

Assume that $h\in \bR$, $\si\in \Re L_2(0,1)$, and let
$(\la_n^2)_{n=0}^\infty$ be the sequence of eigenvalues of the operator
$T=T_{\si,h}$. Denote by $u_n$ the eigenfunction of $T$ corresponding to
the eigenvalue $\la_n^2$ and normalized by the initial condition
$u_n(0)=\sqrt2$. Then $u_n = (I+K_\si) v_n$, where
$v_n(x):=\sqrt2\cos(\la_n x)$. Put
\[
    \al_n:=\|u_n\|_{L_2(0,1)}^{-2};
\]
then the asymptotics of $\la_n$ and properties of the transformation
operator $I+K_\si$ imply that
\[
    \al_n = 1 + \be_n,
\]
where the sequence $(\be_n)_{n\in\bZ_+}$ belongs to
$\ell_2(\bZ_+)$~\cite{HMinv}. Using the Parseval identity
\[
    \slim_{N\to\infty}\sum_{n=0}^N \al_n (\,\cdot\,,u_n) u_n = I
\]
($\slim$ denoting the limit in the strong operator topology in $L_2(0,1)$),
replacing $u_n$ by $(I+K_\si)v_n$, and recalling
relation~\eqref{eq:3.fact}, we conclude that
\[
    I+F_\phi = \slim_{N\to\infty} \sum_{n=0}^N \al_n (\,\cdot\,,v_n)v_n.
\]
Straightforward calculations show that the function $\phi$
of~\eqref{eq:3.phi} determining the kernel~$f_\phi$ of the
operator~$F_\phi$ is given by the series
\begin{equation}\label{eq:3.phi-sum}
    \phi(x) = \sum_{n=0}^\infty [\al_n \cos(\la_nx) - \cos(\pi nx)] +
            \frac12,
\end{equation}
the equality being understood in the $L_2(0,2)$-sense.

\begin{remark}\label{rem:3.coinc}
Assume that $\si_0\in L_2(0,\tfrac12)$ and that $\si$ is an arbitrary
extension of $\si_0$ to a function from $L_2(0,1)$. It is easily seen that
the restriction of the kernel $l_\si$ onto the square
$(0,\tfrac12)\times(0,\tfrac12)$ coincides with the kernel $l_{0,\si_0}$ of
the transformation operator related to the function $\si_0$, see
Section~\ref{sec:form}. Thus the function $\phi_0$ given by
formula~\eqref{eq:2.phi0} verifies the equality
\[
    \phi_0(2x) = -\frac12\si(x) + \int_0^{x} l^2_\si(x,t)\,dt, \quad
        x\in(0,\tfrac12);
\]
in particular, $\phi_0$ is the restriction to $(0,1)$ of the function
$\phi=\phi_\si$ of~\eqref{eq:3.phi}, which implies the inclusion $\phi_0\in
L_2(0,1)$.
\end{remark}

\section{Proof of Theorem~\ref{thm:main}}\label{sec:main}

We start with establishing several lemmata.

\begin{lemma}\label{lem:4.series}
Assume that $(\la_n)_{n\in\bZ_+}\in \fL$ and
$(\be_n)_{n\in\bZ_+}\in\ell_2$. Then the series
\[
    \sum_{n=0}^\infty \be_n \cos(\la_n x)
        \qquad \text{and} \qquad \sum_{n=0}^\infty \be_n \sin(\la_n x)
\]
converge in~$L_2(0,2)$.
\end{lemma}

\begin{proof}
We observe that the systems $\{\cos(\la_nx)\}_{n\in\bZ_+}$ and
$\{\sin(\la_nx)\}_{n\in\bN}$ form Riesz bases of $L_2(0,1)$~\cite{HV}.
Convergence of both series in $L_2(0,1)$ now follows from the definition of
a Riesz basis, see~\cite[Ch.~VI]{GKnsa}. Since
\[
    \be_n\cos[\la_n (x+1)] = \be_n \cos \la_n \cos(\la_n x)
            - \be_n \sin\la_n \sin (\la_nx)
\]
and the sequences $(\be_n\cos\la_n)$ and $(\be_n\sin\la_n)$ belong to
$\ell_2(\bZ_+)$, the series
\[
    \sum_{n=0}^\infty \be_n\cos[\la_n (x+1)]=
         \sum_{n=0}^\infty\be_n \cos \la_n \cos(\la_n x)
        -\sum_{n=0}^\infty\be_n \sin \la_n \sin(\la_n x)
\]
also converges in~$L_2(0,1)$. Henceforth the series
 $\sum_{n=0}^\infty \be_n \cos(\la_n x)$
converges in $L_2(0,2)$. The second series is treated in the same manner,
and the lemma is proved.
\end{proof}

\begin{lemma}\label{lem:4.trig}
For any real numbers $a$ and $b$ the following inequality holds:
\[                      
    |\cos(a+b) - \cos a + b \sin a|\le \frac{b^2}{\sqrt3}.
\]
\end{lemma}

\begin{proof}
Using the standard trigonometric identities, we find that
\begin{align*}          
    |\cos(a+b) - \cos a + b \sin a| &= |(\cos b-1)\cos a + (b-\sin b)\sin a|
        \\
        &\le \sqrt{(\cos b - 1)^2 + (b-\sin b)^2}\\
        &=  \sqrt{2(1-\cos b) + b (b-2\sin b)}.
\end{align*}
Taking into account the inequalities $1-\cos b\le b^2/2$ and
    $|b-\sin b| \le |b|^3/6$ , holding for all real~$b$,
we conclude that
\[                      
    2(1-\cos b) + b (b-2\sin b) \le 2b(b-\sin b) \le b^4/3,
\]
and the result follows.
\end{proof}

\begin{lemma}\label{lem:4.psiLa}
Assume that $(\la_n)\in \fL$ and $(\nu_n)\in\ell_\infty(\bZ_+)$;
then the series
\begin{equation}\label{eq:4.psiLa}      
    \sum_{n=0}^\infty \nu_n[\cos(\la_n x)-\cos (\pi nx)] =:\theta(x)
\end{equation}
converges in $L_2(0,2)$. Moreover, putting $\mu_n:=\la_n-\pi n$
and denoting $\gamma:=\|(\mu_n)\|_{\ell_2}$,
$\delta:=\|(\nu_n)\|_{\ell_\infty}$ we have
\[                      
    \|\theta\|_{L_2(0,1)}\le \delta\Bigl(\frac{\ga^2}{\sqrt{15}} +
        \frac{\ga}{\sqrt2}\Bigr).
\]
\end{lemma}

\begin{proof}
By the definition of the set $\fL$, the sequence $(\mu_n)_{\bZ_+}$ falls
into $\ell_2$. Denote
\[                      
    \ti \theta(x) := -x\sum_{n=0}^\infty \mu_n \nu_n \sin (\pi nx);
\]
then $\ti\theta\in L_2(0,2)$ by Lemma~\ref{lem:4.series} and,
since the set $\{\sqrt2\sin(\pi nx)\}_{n\in\bN}$ is an orthonormal
basis of $L_2(0,1)$,
\[                      
    \|\ti\theta\|_{L_2(0,1)}\le
        \Bigl(\frac12 \sum_{n=0}^\infty |\mu_n\nu_n|^2\Bigr)^{1/2}
        \le \frac{\ga\delta}{\sqrt2}.
\]
Applying Lemma~\ref{lem:4.trig}, we find that
\begin{equation}\label{eq:4.8}
    |\cos(\la_n x) - \cos(\pi nx) + \mu_n x \sin(\pi nx)|
        \le \frac{|\mu_nx|^2}{\sqrt3}, \qquad x\in[0,2],
\end{equation}
so that the series
\[
    \sum_{n=0}^\infty \nu_n\bigr[ \cos(\la_n x) - \cos(\pi nx)
         + \mu_n x \sin(\pi nx)\bigl]
\]
converges uniformly and absolutely on $[0,2]$. This proves that
series~\eqref{eq:4.psiLa} converges in $L_2(0,2)$; we denote its
sum by $\theta$.

Inequality~\eqref{eq:4.8} yields the estimate
\[
    |\theta(x) - \ti\theta(x)|
        = \Bigl| \sum_{n=0}^\infty \nu_n \bigr[ \cos(\la_n x) - \cos(\pi nx)
             + \mu_n x \sin(\pi nx)\bigl]\Bigr|
         \le \frac{\ga^2\delta}{\sqrt3}x^2,
\]
so that $\|\theta - \ti\theta\|_{L_2(0,1)}\le
\ga^2\delta/{\sqrt{15}}$ and
\[
    \|\theta\|_{L_2(0,1)} \le \|\theta - \ti\theta\|_{L_2(0,1)}
        + \|\ti\theta\|_{L_2(0,1)} \le \frac{\ga^2\delta}{\sqrt{15}} +
        \frac{\ga\delta}{\sqrt2}
\]
as claimed. The proof is complete.
\end{proof}

\begin{lemma}\label{lem:4.Uinv}
Assume that $(\la_n)\in\fL$ and that $(\al_n)$ is a sequence of real
numbers such that $c<\al_n<C$ for some positive constants~$c,C$ and all
$n\in\bZ_+$. Define an operator $U:\,L_2(0,1) \to L_2(0,1)$ by the equality
\begin{equation}\label{eq:4.U}
    U:= \slim_{N\to\infty}\sum_{n=0}^N \al_n (\,\cdot\,,v_n)v_n,
\end{equation}
where $\slim$ stands for the limit in the strong operator topology
of~$L_2(0,1)$ and $v_n(x):=\sqrt2\cos(\la_nx)$. Then the operator $U$ is
invertible and
 \(
    (U^{-1} v_j,v_k) = \al_k^{-1} \delta_{jk}
 \)
for all $j,k\in\bZ_+$, where $\delta_{jk}$ is the Kronecker delta.
\end{lemma}

\begin{proof}
Observe that in view of the relation $\la_n-\pi n\to0$ as $n\to\infty$ the
system $\{\cos(\la_nx)\}_{n\in\bZ_+}$ is a Riesz basis of
$L_2(0,1)$~\cite{HV}, which yields convergence of the sum in~\eqref{eq:4.U}
in the strong operator topology. It is easily seen that
\[
    U^{-1} = \slim_{N\to\infty} \sum_{n=1}^N \al_n^{-1}
        (\,\cdot\,, \hat v_n)\hat v_n,
\]
where $(\hat v_n)$ is a basis biorthogonal to $(v_n)$,
see~\cite[Ch.~VI]{GKnsa}. Therefore
\[
    (U^{-1}v_j,v_k) = \slim_{N\to\infty}
        \sum_{l=1}^N \al_l^{-1} ( v_j , \hat v_l)(\hat v_l, v_k)
    = \al_l^{-1} \delta_{jl}\delta_{lk} = \al_k^{-1} \delta_{jk},
\]
and the proof is complete.
\end{proof}

\begin{proof}[Proof of Theorem~\ref{thm:main}.]
(i) \emph{Necessity.} Assume that $\si_0\in \Re L_2(0,\tfrac12)$,
$(\la_n)=:\La\in\fL$, and let there exist an extension~$\si\in L_2(0,1)$
of~$\si_0$ to~$(0,1)$ and a number $h\in\bR$ such that the spectrum of the
corresponding operator $T_{\si,h}$ coincides with the set $\{\la_n^2\}$.
Then according to Remark~\ref{rem:3.coinc} and formula~\eqref{eq:3.phi-sum}
we have
\begin{equation}\label{eq:4.1}
    \phi_0(x) = \sum_{n=0}^\infty [\al_n\cos(\la_nx) - \cos(\pi nx)]
        + \tfrac12, \qquad x\in(0,1)
\end{equation}
for some $\al_n>0$ such that $\be_n:=\al_n-1$ form an $\ell_2$-sequence, so
that $\phi_0\in \Pi_\La$ as required.

\emph{Sufficiency.} Assume that $\si_0\in \Re L_2(0,\tfrac12)$,
$(\la_n)=:\La\in\fL$, and that the function $\phi_0$ belongs to $\Pi_\La$.
According to the definition of the set $\Pi_\La$, the function~$\phi_0$ has
the form~\eqref{eq:4.1} for given $\la_n$ and some sequence $(\al_n)$ of
positive numbers, for which $\be_n:=\al_n-1$ form an $\ell_2$-sequence.
Writing $\al_n\cos(\la_nx)-\cos(\pi nx)$ as
\[
    \be_n\cos(\la_nx) + [\cos(\la_nx)-\cos(\pi nx)]
\]
and applying Lemmata~\ref{lem:4.series} and \ref{lem:4.psiLa}, we conclude
that the series on the right-hand side of~\eqref{eq:4.1} converges in
$L_2(0,2)$ to some function $\phi$; clearly, $\phi(x) = \phi_0(x)$ a.e. on
$(0,1)$.

Consider now the operator $I+F_\phi$ corresponding to the function $\phi$
constructed above. Using the definition of $F_\phi$, one easily shows that
\[
    I+ F_\phi = \slim_{N\to\infty}
        \sum_{n=0}^N \al_n (\,\cdot,v_{n})v_{n}
\]
with $v_n(x) = \sqrt2\cos(\la_nx)$. Since the set $\{v_{n}\}_{n\in\bZ_+}$
forms a Riesz basis of $L_2(0,1)$ and the numbers $\al_n$ are uniformly
positive and uniformly bounded, the operator $I+F_\phi$ is bounded and
(uniformly) positive. It follows that the operator $I+F_\phi$ is
factorizable~\cite[Ch.~4]{GKvolt}, so that  $\phi\in\Phi$ and by
Theorem~\ref{thm:3.conn} there exists $\si\in L_2(0,1)$ such that
\begin{equation}\label{eq:sifound}
    \phi(2x) = -\frac12\si(x) + \int_0^x l_\si^2(x,t)\,dt,
        \quad x\in(0,1).
\end{equation}
Here $l_\si$ is the kernel of the transformation operator corresponding to
$\ti T_0$ and $\ti T_{\si}$.

We prove next that $\si=\si_0$ on $(0,\tfrac12)$. Denote by $\ti\si_0$ the
restriction of $\si$ to $(0,\tfrac12)$; then by Remark~\ref{rem:3.coinc}
\[
    l_\si(x,t) = l_{0,\ti\si_0}(x,t), \qquad 0<t<x<\tfrac12,
\]
where $l_{0,\ti\si_0}$ is the kernel of the transformation operator on the
interval $(0,\tfrac12)$ constructed for the function $\ti\si_0$. On the
other hand, $\phi_0$ satisfies the relation
\[
    \phi_0(2x) = -\frac12\si_0(x) + \int_0^x l_{0,\si_0}^2(x,t)\,dt,
        \quad x\in(0,\tfrac12),
\]
and since the mapping~\eqref{eq:3.phi} considered on $(0,\tfrac12)$ is
bijective, we conclude that $\ti\si_0 = \si_0$, i.e., that $\si$ is an
extension of $\si_0$.

It remains to show that there is $h\in\bR$ such that the spectrum of
$T_{\si,h}$ coincides with the set $\{\la_n^2\}_{n\in\bZ_+}$. Put
\[
    w_j:= (I+K_\si)v_{j}, \quad j\in \bZ_+;
\]
then by virtue of Lemma~\ref{lem:4.Uinv} we find that
\[
    (w_j,w_k) =
        \bigl((I+K_\si^*)(I+K_\si) v_{j},v_{k}\bigr)
        =\bigl((I+F_\phi)^{-1}v_{j},v_{k}\bigr) = \al_k^{-1}\delta_{jk}.
\]
In particular, the vectors $w_j$ are orthogonal; since
$\{v_{k}\}_{n\in\bZ_+}$ is a Riesz basis of $L_2(0,1)$ and the
operator $I+K_\si$ is a homeomorphism of $L_2(0,1)$, the set
$\{w_j\}_{n\in\bZ_+}$ is an orthogonal Riesz basis of $L_2(0,1)$.

It remains to prove that there is $h\in\bR$ such that $w_j$ are
eigenfunctions of the operator $T_{\si,h}$, i.e., that $w_j^{[1]}=hw_j(1)$
for all $j\in\bZ_+$. We observe that $w_j$ satisfy the following Lagrange
identity:
\begin{equation}\label{eq:4.16}
    0 = (\ti T_{\si} w_j,w_k) - (w_j, \ti T_{\si} w_k)
        =  - w_j^{[1]}(1)\ov{w_k(1)} + w_j(1)\ov{w^{[1]}_k(1)}.
\end{equation}
If $w_j(1)=0$ for some $j\in\bZ_+$, then $w_j^{[1]}(1)\ne0$ by uniqueness
theorem for the equation $\ell_\si(y) = \la_j^2y$~\cite{SS} and the above
relations imply that $w_k(1)=0$ for all $k\in\bZ_+$. Then~$\la_k^2$ would
be eigenvalues of the Sturm--Liouville operator $T_{\si,\infty}$. This is
impossible since the eigenvalues $\nu_k^2$ of $T_{\si,\infty}$ obey the
asymptotics $\nu_k = \pi (n-\tfrac12) + \mathrm{o}(1)$ as
$k\to\infty$~\cite{HMinv,SS1}, and this asymptotics is completely different
from that of $\la_k$. Henceforth $w_k(1)\ne0$ for all $k\in\bZ_+$, so that
by~\eqref{eq:4.16} there exists $h\in\bR$ such that $w_k^{[1]}(1)/w_k(1) =
h$, i.e., such that the set $\{\la_k^2\}$ are eigenvalues of the operator
$T_{\si,h}$ and $w_k$ are the corresponding eigenfunctions. The operator
$T_{\si,h}$ has no other eigenvalues since the set $\{w_j\}_{j\in\bZ_+}$ of
the eigenfunctions is already complete in $L_2(0,1)$.

To prove (ii) we observe that the spectral data $(\si_0,\La)$ uniquely
determine the transformation operator $K$, so that $\si'=q$ is unique. The
above reasonings show that the number $h$ in the boundary condition is then
identified as $w_0^{[1]}(1)/w_0(1) = h$, where $w_0 = (I+K)v_0$ and
$v_{0}=\sqrt{2}\cos(\la_0x)$. Thus the extension $q$ of $q_0$ and $h\in\bR$
are unique, and the proof is complete.
\end{proof}

Next we give an example showing that the inclusion $\phi_0\in\Pi_\La$ need
not hold in general, so that all the hypotheses of Theorem~\ref{thm:main}
are essential for solvability of the half-inverse problem.

\begin{example}\label{ex:4.ex-solv}
Set $\la_n:=\pi n$, $n\in\bZ_+$, $\La:=(\la_n)_{n\in\bZ_+}$, and
\begin{equation}\label{eq:4.siga0}
    \si_{\ga,0}(x):=\frac{2\ga}{1-\ga x}-\ga,
        \qquad x\in [0,\tfrac12].
\end{equation}
Observe that $\si_{\ga,0}$ belongs to $L_2(0,\tfrac12)$ if
$\ga<2$. We shall show that the corresponding function
$\phi_{\ga,0}:=\phi_0(\cdot,\si_{\ga,0})$ of~\eqref{eq:2.phi0}
equals $-\gamma/2$; therefore $\phi_{\ga,0}\in\Pi_\La$ if and only
if $\ga<1$, so that for $\ga<2$ the half-inverse spectral problem
with mixed data $\{\si_{\ga,0},\La\}$ is soluble if and only if
$\ga<1$.

By a straightforward verification one sees that a solution
$y_0(\cdot,\la)$ to the equation
\[
    -(y' - \si_{\ga,0}y)' - \si_{\ga,0}y' = \la^2 y
\]
satisfying the initial conditions
    $y(0)=1$, $(y' - \si_{\ga,0}y)(0)=0$, is
\[
    y_0(x,\la)= \cos(\la x) + \frac{\ga}{\la}\frac{\sin(\la x)}{1-\ga x}.
\]
The kernel $l_{0,\si_{\ga,0}}(x,t)$ of the transformation operator
$I+L_{0,\si_{\ga,0}}$ must satisfy the following identity for all
$\la\in\bC$ and all $x\in[0,\tfrac12]$:
\[
    \cos(\la x) + \frac{\ga}{\la}\frac{\sin(\la x)}{1-\ga x}
        + \int_0^x l_{0,\si_{\ga,0}}(x,t)
        \Bigl[\cos(\la t)
            + \frac{\ga}{\la}\frac{\sin(\la t)}{1-\ga t}\Bigr]\,ds
        = \cos(\la x)
\]
Observe that
\[
    \left(\frac{\ga}{\la}\frac{\sin(\la x)}{1-\ga x}\right)' =
    \frac{\ga}{1-\ga x}\Bigl[\cos(\la x)
            + \frac{\ga}{\la}\frac{\sin(\la x)}{1-\ga x}\Bigr],
\]
which suggests that
\[
    l_{0,\si_{\ga,0}}(x,t)=-\frac{\ga}{1-\ga t}.
\]
Now the function $\phi_{\ga,0}$ is found to be
\begin{align*}
    \phi_{\ga,0}(2x) &=  - \frac12\si_{\ga,0}(x)
        + \int_0^{x} l^{\,2}_{0,\si_{\ga,0}}(x,t)\,dt =
        -\frac{\ga}{1-\ga x}+\frac{\ga}2
        + \int_0^x \frac{\ga^2}{(1-\ga t)^2}\,dt\\
    &= -\frac{\ga}{1-\ga x}+\frac{\ga}2 +\frac{\ga}{1-\ga t}\Bigr\vert_{t=0}^{t=x}
        =  -\frac{\ga}2, \qquad x\in(0,\tfrac12),
\end{align*}
so that indeed $\phi_{\ga,0}\in \Pi_\La$ if and only if $\ga<1$.

For $\ga<1$, the half inverse problem can be solved explicitly.
Indeed, according to the reconstruction algorithm of
Theorem~\ref{thm:main} the extension $\phi_\ga$ of the function
$\phi_{\ga,0}$ to the interval $(0,2)$ is given by the same
formula, i.e., $\phi_{\ga}\equiv-\ga/2$ on $(0,2)$, and the
integral operator $F_{\phi_\ga}$ corresponding to $\phi_\ga$ has
kernel
\[
    f_{\phi_\ga}(x,t) \equiv -\ga, \quad x,t\in [0,1].
\]
The operator $I+ F_{\phi_\ga}$ is selfadjoint; moreover, for
$\ga<1$ it is uniformly positive, and thus $\phi_\ga\in\Phi$ for
such $\ga$.

A simple observation suggests that a solution $k$ to the
Gelfand--Levitan--Marchenko equation
\[
    k(x,t) + f_{\phi_\ga}(x,t) + \int_0^x k(x,s)f_{\phi_\ga}(s,t)\,ds = 0,
        \quad 0<t<x<1,
\]
must have the form $k(x,t)=a(x)$; after a straightforward
calculation we conclude that
\[
    k(x,t) = \frac{\ga}{1-\ga x},   \quad 0<t<x<1.
\]
By~\cite{HMtr,HMinv} $k$ is a kernel of the transformation operator $K_\si$
with
\[
    \si(x) = 2 k(x,x) + 2 \phi(0) = \si_\ga(x),
\]
where the function $\si_\ga$ is the extension of $\si_{\ga,0}$ to the
interval~$(0,1)$ by formula~\eqref{eq:4.siga0}. The constant $h$ in the
boundary conditions is identified as $w_0^{[1]}(1)/w_0(1)$, with
\[
    w_0(x):= (I+K_\si)1 = 1 + \frac{\ga x}{1-\ga x}
            = \frac{1}{1-\ga x}.
\]
In particular, $w_0(1)=1/(1-\ga)$, $w'_0(1)=\ga/(1-\ga)^2$,
$\si(1)=(\ga+\ga^2)/(1-\ga)$, so that $h_\ga=-\ga^2/(1-\ga)$. Now
it can be directly checked that the spectrum of the operator
$T_{\si_\ga,h_\ga}$ coincides with the set
$\{\pi^2n^2\}_{n\in\bZ_+}$, the corresponding eigenfunctions being
$w_0$ above and
\[
    w_n(x) = (I+K_\si)\cos(\pi n x)=
        \cos(\pi n x)+ \frac{\ga}{\pi n}\frac{\sin(\pi n x)}{1-\ga x},
        \quad n\in\bN.
\]
Thus the half-inverse problem is solved.
\end{example}


\section{Proof of Theorems~\ref{thm:local}
    and~\ref{thm:regular}}\label{sec:prf2.2}

We proved in Theorem~\ref{thm:main} that the half-inverse problem
is not soluble for any mixed data $\{\si_0,\La\}\in\Re
L_2(0,\tfrac12)\times\fL$. However, the set of those mixed data,
for which a solution exists, can easily be shown to be open in
$\Re L_2(0,\tfrac12)\times\fL$, the topology in $\fL$ being
inherited from that of $\ell_2(\bZ_+)$ through identification of
$\La=(\la_n)\in\fL$ and $(\la_n-\pi n)\in \ell_2(\bZ_+)$. This
follows from the fact that both the mapping $\si_0\mapsto \phi_0$
of~\eqref{eq:2.phi0} and the mapping $\La\mapsto \psi_\La$ induced
by~\eqref{eq:2.psiLa} are continuous. In particular, there exists
a neighbourhood in $\Re L_2(0,\tfrac12) \times \fL$ of the
``unperturbed'' mixed data $\si_0\equiv 0$, $\La = (\pi
n)_{n\in\bZ_+}$, in which the half-inverse problem is soluble. The
aim of Theorem~\ref{thm:local} is to make this observation
quantitative, and to this end it suffices to estimate the norms of
the functions $\phi_0$ and $\psi_\La$ in terms of $\si_0$ and
$\La$. As we mentioned in Section~\ref{sec:form}, it is easier to
restrict ourselves to the case $\si_0\in \Re W^1_2(0,\tfrac12)$,
i.e., to the case $q_0\in \Re L_2(0,\tfrac12)$, and use the
corresponding norms.

We start with the following auxiliary lemma.

\begin{lemma}\label{lem:5.phi0}
Assume that $q_0\in L_2(0,\tfrac12)$ is such that $\|q_0\|_{L_2(0,1/2)} \le
\tfrac12$. Put $\si_0(x):=\int_0^x q_0(t)\,dt$ for $x\in(0,\tfrac12)$ and
construct the function $\phi_0$ on $(0,1)$ via~\eqref{eq:2.phi0}. Then
$\|\phi_0\|_{L_2(0,1)}\le \tfrac14$.
\end{lemma}

\begin{proof}
Formula~\eqref{eq:2.phi0} suggests that to bound the norm of the function
$\phi_0$ we only need to estimate the kernel $l_{0,\si_0}$ of the
transformation operator $L_{0,\si_0}$ constructed as explained in
Section~\ref{sec:form}. In the case where $q_0\in L_2$ the differential
expression $\ell_{\si_0}(u):= -(u'-\si_0 u)' - \si u'$ coincides with $-u''
+ q_0u$ and it is well known (see~\cite[Ch.~1.2]{Ma}) that the kernel
$l_0=l_{0,\si_0}$ is then a unique solution to the hyperbolic equation
\[
    l''_{xx}(x,t) = l''_{tt}(x,t) - q_0(t)l(x,t)
\]
satisfying the boundary conditions
\[
    l(x,x) = -\frac12\int_0^x q_0(t)\,dt= -\frac12\si_0(x), \qquad l'_t(x,t)\vert_{t=0} = 0.
\]

Introduce new variables $u=x+t$ and $v=x-t$ and put
 \(
    a(u,v):= l\bigl(\tfrac{u+v}2,\tfrac{u-v}2 \bigr);
 \)
then standard reasonings (see, e.g., \cite[Ch.~1.2]{Ma}) reduce the above
boundary value problem to the integral equation
\begin{equation}\label{eq:5.inteq}
\begin{aligned}
    a(u,v) &= -\frac14\int_{v}^{u}d\al\int_0^{v}
            q_0\Bigl(\frac{\al-\be}2\Bigr)  a(\al,\be)\,d\be\\
           &\quad -\frac12\int_{0}^{v}d\al\int_0^{\al}
            q_0\Bigl(\frac{\al-\be}2\Bigr)  a(\al,\be)\,d\be\\
    &\quad -\frac12\biggl[\int_0^{u/2}q_0(\al)\,d\al
                + \int_0^{v/2}q_0(\be)\,d\be\biggr].
\end{aligned}
\end{equation}
This integral equation possesses a unique solution $a$, which is continuous
on the set $0\le v \le u \le 1$. We put
\[
    b(s) := \max_{0\le v\le u\le s} |a(u,v)|,
        \quad s\in[0,1],
\]
and observe that for any fixed $s\in[0,1]$ the above maximum is assumed at
some point $(u_0,v_0)$ satisfying the relation $0\le v_0 \le u_0 \le s$.
Integral equation~\eqref{eq:5.inteq} then implies that
\begin{align*}
    b(s) &\le \frac{b(s)}2 \int_{0}^{s}\,d\al\int_0^{\al}
            \Bigl|q_0\Bigl(\frac{\al-\be}{2}\Bigr)\Bigr|\,d\be
            + \int_0^{s/2}  |q_0(\al)|\,d\al \\
         &\le [s b(s) + 1 ] \int_0^{s/2}|q_0(\al)|\,d\al
             \le [b(s) + 1] \int_0^{1/2}|q_0(\al)|\,d\al.
\end{align*}
Taking into account the inequality
\[
    \int_0^{1/2}|q_0(\al)|\,d\al \le \frac1{\sqrt{2}}\|q_0\|_{L_2(0,1/2)}
        \le \frac1{2\sqrt2},
\]
we conclude that
\[
    b(s) \le \frac1{2\sqrt2}\Bigl(1-\frac1{2\sqrt2}\Bigr)^{-1}
        =\frac1{2\sqrt{2}-1} < \frac1{\sqrt3}.
\]
Therefore for the solution $l_0$ of integral equation~\eqref{eq:5.inteq} we
find that
\[
    \phi_1(x):=\int_0^{x/2} l_0^2(x/2,t)\,dt \le \frac{x}6,
\]
so that
 \(
    \|\phi_1\|_{L_2(0,1)}\le \frac1{6\sqrt3}.
 \)
Using the estimate
\[
    |\si_0(x/2)| \le \int_0^{x/2}|q_0(\al)|\,d\al
        \le \sqrt{\frac{x}2}\|q_0\|_{L_2(0,1/2)}\le
        \frac{\sqrt{x}}{2\sqrt2},
\]
we arrive at
\[
    \int_0^1 |\si_0(x/2)|^2\,dx \le \frac18\int_0^1 x\,dx = \frac1{16}.
\]
Finally, \eqref{eq:2.phi0} yields the required bound
\[
    \|\phi_0\|_{L_2(0,1)} \le \frac12\|\si_0(\cdot/2)\|_{L_2(0,1)}
        + \|\phi_1\|_{L_2(0,1)} \le \frac18 + \frac1{6\sqrt3} \le \frac14,
\]
and the proof is complete.
\end{proof}

\begin{proof}[Proof of Theorem~\ref{thm:local}]
We have to show that under the hypotheses of the theorem the
function~$\phi_0$ of~\eqref{eq:2.phi0} belongs to~$\Pi_\La$. As explained
in Section~\ref{sec:form}, the function $\phi_0$ can be represented in the
form
\[
    \phi_0(x) = \sum_{n=0}^\infty [\al_n \cos(\la_nx)-\cos(\pi nx)] +\tfrac12,
\]
where the real numbers $\al_n$ are such that the sequence
$(\be_n)_{n\in\bZ_+}$ with $\be_n:=\al_n-1$ belongs to $\ell_2$.
Thus the only thing to be proved is that all $\al_n$ are positive.

Set
\[
 \theta(x) :=\sum_{n=0}^\infty \alpha_n[\cos (\pi nx) -\cos(\la_nx)],
    \qquad x\in(0,1);
\]
by Lemma~\ref{lem:4.psiLa} the above series converges in
$L_2(0,1)$, so that $\theta\in L_2(0,1)$. Then
\[
    \phi_0(x) + \theta(x) = \sum_{n=0}^\infty (\alpha_n-1)\cos (\pi nx)
    +\tfrac12 = \sum_{n=1}^\infty (\alpha_n-1)\cos (\pi nx)
    +\al_0 - \tfrac12,
\]
so that
\[
    \|\phi_0 + \theta\|^2_{L_2(0,1)} =
        \frac12 \sum_{n=1}^\infty (\alpha_n-1)^2
            +(\alpha_0-\tfrac12)^2.
\]
On the other hand, $\|\phi_0\|_{L_2(0,1)} \le \tfrac14$ by
Lemma~\ref{lem:5.phi0}, and the norm of the function $\theta$ can be
estimated by virtue of Lemma~\ref{lem:4.psiLa} as
\[
    \|\theta\|_{L_2(0,1)} \le \delta\Bigl(\frac{\ga^2}{\sqrt{15}} +
                \frac{\ga}{\sqrt2}\Bigr)
        \le  \frac{\delta}{4}\Bigl(\frac{1}{\sqrt 2}+ \frac{1}{4\sqrt 15}\Bigr)
        \le \frac{\delta}{4}\Bigl(\frac{11}{15}+  \frac{1}{15}\Bigr)
        = \frac{\delta}{5},
\]
where we have put $\ga:=\|(\mu_n)\|_{\ell_2}$,
$\delta:=\sup_{n\in\bZ+}|\al_n|$, and used the fact that $\ga\le\tfrac14$
by the assumption of the theorem. Taking into account the above relations,
we conclude that
\begin{equation}\label{eq:5.norm}
    \frac12 \sum_{n=1}^\infty (\alpha_n-1)^2  +(\alpha_0 -\tfrac12)^2
        \le \Bigl(\frac{1}{4}+ \frac{\delta}{5}\Bigr)^2,
\end{equation}
which yields the estimates
\begin{equation}\label{eq:5.al}
    |\alpha_0 -\tfrac{1}{2}| \le \frac{1}{4}+
     \frac{\delta}{5}, \qquad
    |\alpha_n -1| \le \sqrt 2 \Bigl(\frac{1}{4}+ \frac{\delta}{5}\Bigr),
            \qquad n\in \bN.
\end{equation}
These inequalities imply that
\[
    \delta \le 1+\sqrt2 \Bigl(\frac{1}{4}+
        \frac{\delta}{5}\Bigr) \le \frac{14}{10}+ \frac{3\delta}{10},
\]
i.e., that $\delta \le 2$. Returning now to
inequalities~\eqref{eq:5.al}, we get
\[
    |\alpha_n -1| \le \frac{2\sqrt 2}{3} <1,\qquad n\in\bN,
\]
so that $\alpha_n >0$ for all $n\in \bN$.

Assume that $\al_0<0$. This is only possible when $\delta>1$ (as otherwise
$|\al_0-\tfrac12|<\tfrac12$ and $\al_0>0$), and then
$\delta=\sup_{n\in\bN}\al_n$. Therefore relation~\eqref{eq:5.norm} yields
the inequality
\[
     \tfrac12 (\delta-1)^2+ (\alpha_0 -\tfrac12)^2 \le \left(\frac{1}{4}+
            \frac{\delta}{5}\right)^2,
\]
so that, on account of the assumption $\al_0<0$, we must have
\[
    P(\delta):= \frac12 (\delta-1)^2 + \frac14 -
        \left(\frac{1}{4}+ \frac{\delta}{5}\right)^2 \le 0.
\]
However, the last inequality never holds for real $\delta$ as the
discriminant of the polynomial~$P$ is negative. The contradiction derived
shows that the assumption $\al_0<0$ was wrong, and the theorem is proved.
\end{proof}

Before proceeding with the proof of Theorem~\ref{thm:regular}, we refine
the statements of Lemmata~\ref{lem:4.psiLa} and \ref{lem:4.series} for the
case where $\La\in \fL_1$.

\begin{lemma}\label{lem:5.psiLasmooth}
Assume that $\La\in\fL_1$; then the function~$\psi_\La$
of~\eqref{eq:2.psiLa} belongs to~$W^1_2(0,2)$.
\end{lemma}

\begin{proof}
We recall that $\psi_\La$ is given by the series
\[
    \sum_{n=0}^\infty [\cos(\la_n x)-\cos (\pi nx)],
\]
which converges in the topology of~$L_2(0,2)$. We write
\begin{align*}
    \cos(\la_nx) - \cos(\pi nx)
        &= \Bigl[\cos(\mu_n x) -1 + \frac{\mu_n^2x^2}{2}\Bigr]\cos(\pi nx)
            - [\sin(\mu_nx) - \mu_nx]\sin(\pi nx)\\
        &\qquad  - \frac{\mu^2_nx^2}{2}\cos(\pi nx) - \mu_nx\sin(\pi nx)\\
        &=: \om_{1,n}(x) - \om_{2,n}(x)
            - \frac{\mu^2_nx^2}{2}\cos(\pi nx) - \mu_nx\sin(\pi nx),
\end{align*}
where as usual $\mu_n$ equals $\la_n-\pi n$.

By definition, $\La\in\fL_1$ means that $\mu_n = c/n + \nu_n/n$
for some $c\in\bR$ and an $\ell_2$-sequence~$(\nu_n)$, and this
representation implies that the series
$\sum_{n\in\bZ_+}\mu_n\sin(\pi nx)$ converges in $L_2(0,2)$ to a
function from $W^1_2(0,2)$. Indeed, we have
\[
    \sum_{n=1}^\infty \frac{\sin(\pi nx)}n =\frac{\pi}2(1-x),
        \quad x\in(0,2),
\]
the equality being understood in the $L_2$-sense, while the series
\[
    \sum_{n=1}^\infty \frac{\nu_n}n \sin(\pi nx)
\]
converges in $W^1_2(0,2)$. Similarly, the series $\sum_{n=0}^\infty \mu_n^2
\cos(\pi nx)$ converges in $W^1_2(0,2)$. Next, we observe that
\[
    \om'_{1,n}(x) = {\mathrm O}(n^{-3}), \qquad
    \om'_{2,n}(x) = {\mathrm O}(n^{-2})
\]
as $n\to\infty$, so that the series
\[
    \sum_{n=0}^\infty \om_{1,n}(x), \qquad
    \sum_{n=0}^\infty \om_{2,n}(x)
\]
converge in $C^1[0,2]$. Summing up, we see that the function
\[
    \psi_\La = \sum_{n=0}^\infty \om_{1,n}
            - \sum_{n=0}^\infty \om_{2,n}
            - \frac{x^2}2 \sum_{n=0}^\infty \mu_n^2 \cos(\pi nx)
            - x \sum_{n=0}^\infty \mu_n\sin(\pi nx)
\]
belongs to $W^1_2(0,2)$, and the lemma is proved.
\end{proof}

\begin{lemma}\label{lem:5.psi1}
Assume that $\La\in\fL_1$ and that a sequence $(c_n)_{n\in\bZ_+}$ from
$\ell_2$ is such that the series
\begin{equation}\label{eq:5.cn-ser}
    \sum_{n\ge0} c_n \cos(\la_nx)
\end{equation}
converges in $L_2(0,1)$ to a function from $W^1_2(0,1)$. Then this series
converges in $L_2(0,2)$ to a function from $W^1_2(0,2)$.
\end{lemma}

\begin{proof}
Arguing as in the proof of Lemma~\ref{lem:5.psiLasmooth}, we can show that
under the assumptions of the lemma the series
\begin{equation}\label{eq:5.cn-dif-ser}
    \sum_{n\ge0} c_n [\cos(\la_nx) - \cos \pi nx]
\end{equation}
converges in $L_2(a,a+1)$ to a function from $W^1_2(a,a+1)$ for any
 real~$a$. In particular, $c_n$ are cosine Fourier coefficients of
a $W^1_2(0,1)$-function and thus by integration by parts are easily shown
to have the form
\[
    c_n = \frac{c}{n} + \frac{\nu_n}{n}
\]
with $c\in\bR$ and a sequence $(\nu_n)$ falling into $\ell_2$.
This implies that the series
 \[
    \sum_{n\ge0} c_n \cos(\pi n x)
 \]
converges in $L_2(a,a+1)$ to a function from $W^1_2(a,a+1)$ for any real
$a$ (cf. the proof of Lemma~\ref{lem:5.psiLasmooth}). Since the same
statement holds for series~\eqref{eq:5.cn-dif-ser}, the required
convergence result for series~\eqref{eq:5.cn-ser} on the interval $(0,2)$
follows.
\end{proof}

\begin{proof}[Proof of Theorem~\ref{thm:regular}]
We have to prove that in the case where $\si_0\in
W_2^1(0,\tfrac12)$, $\La\in \fL_1$, and $\phi_0\in \Pi_\La$ the
function $\si$ constructed in Theorem~\ref{thm:main} and
extending~$\si_0$ belongs to~$W_2^1(0,1)$. In view of
Theorem~\ref{thm:3.smooth} this is equivalent to showing that the
function $\phi$ of~\eqref{eq:3.phi} belongs to $W_2^1(0,2)$.

Observe first that the function~$\phi_0$ of~\eqref{eq:2.phi0} is the
restriction of~$\phi$ to the interval~$(0,1)$ and that $\phi_0$ belongs to
$W^1_2(0,1)$. To prove this, we extend $\si_0$ to some function $\ti \si$
from~$W^1_2(0,1)$ and construct the function~$\ti \phi$ corresponding
to~$\ti\si$ according to~\eqref{eq:3.phi}. By Theorem~\ref{thm:3.smooth}
$\ti\phi$ is in~$W^1_2(0,2)$, and it remains to notice that, in view of
Remark~\ref{rem:3.coinc}, $\phi_0$ is the restriction of $\ti\phi$
onto~$(0,1)$, so that $\phi_0\in W^1_2(0,1)$.

In virtue of the results of Section~\ref{sec:transf} the function~$\phi$
of~\eqref{eq:3.phi} has the form~\eqref{eq:3.phi-sum}, i.e.,
\[
    \phi(s) = \sum_{n=0}^\infty [\al_n \cos(\la_ns) - \cos(\pi ns)]
            + \tfrac12.
\]
The restriction $\phi_0$ of $\phi$ to $(0,1)$ is therefore given in
$L_2(0,1)$ by the same series. The assumption $\phi_0\in \Pi_\La$ now
implies that $\al_n=1+\be_n>0$ with $(\be_n)\in \ell_2$, so that $\phi$ can
be represented as
\begin{align*}
    \phi(s) &=  \sum_{n=0}^\infty[\cos(\la_ns)-\cos(\pi ns)]
        + \sum_{n=0}^\infty \be_n \cos(\la_ns) + \tfrac12\\
        &=:  \psi_\La(s) + \psi_1(s) + \tfrac12.
\end{align*}

Since $\La\in\fL_1$, by virtue of Lemma~\ref{lem:5.psiLasmooth} the
function $\psi_\La$ belongs to $W_2^1(0,2)$. Restricting the above equality
to $(0,1)$ and recalling that $\phi_0\in W^1_2(0,1)$, we see that
$\psi_1\in W^1_2(0,1)$. Applying now Lemma~\ref{lem:5.psi1} to the series
$\sum_{n\ge0}\be_n\cos(\la_ns)$, we conclude that this series converges in
$L_2(0,2)$ to a function from $W^1_2(0,2)$. Thus $\psi_1\in W^1_2(0,2)$,
and hence $\phi\in W^1_2(0,2)$ as required. The theorem is proved.
\end{proof}

\appendix
\section{Proof of Theorem~\ref{thm:3.smooth}}\label{sec:smooth}

We denote by $G_2$ the subspace of all functions $k\in
L_2\bigl((0,1)\times(0,1)\bigr)$, for which the mappings
\begin{equation}\label{eq:A.xtoL2}
    x \mapsto k(x,\,\cdot\,) \in L_2(0,1), \qquad
    t \mapsto k(\,\cdot\,,t) \in L_2(0,1)
\end{equation}
belong to $C\bigl([0,1],L_2(0,1)\bigr)$ (i.e., for which these mappings
coincide a.e.~with continuous ones). We also denote by $G^1_2$ the subspace
of~$G_2$ consisting of those $k$, for which mappings~\eqref{eq:A.xtoL2}
belong to $C^1\bigl([0,1],L_2(0,1)\bigr)$. The spaces of integral operators
with kernels from $G_2$ (respectively, from $G_2^1$) will be denoted
by~$\fG_2$ (respectively, by~$\fG^1_2$).

\begin{lemma}\label{lem:A.Fphi}
Assume that $\phi\in W^1_2(0,2)$; then the operator $F_\phi$ with kernel
$f_\phi(x,t):=\phi(x+t)+ \phi(|x-t|)$ belongs to $\fG_2^1$.
\end{lemma}

\begin{proof}
Since the function $k(x,t):=\phi'(x+t)$ obviously belongs to
$G_2$, we use the equalities
\[
    \phi(x+t) = \phi(t) + \int_0^x \phi'(t+ \xi)\,d\xi
              = \phi(x) + \int_0^t \phi'(x+ \xi)\,d\xi
\]
to justify the inclusion $\phi(x+t) \in G^1_2$. In the same manner we show
that $\phi(|x-t|)| \in G_2^1$, which by definition yields $F_\phi\in
\fG_2^1$.
\end{proof}

\begin{lemma}\label{lem:A.invert}
Assume that $R\in\fG^1_2$ and that the operator $I+R$ is invertible. Then
the integral operator $\hat R:= (I+R)^{-1} - I$ belongs to $\fG^1_2$.
\end{lemma}

\begin{proof}
We use a kind of the bootstrap method based on the formula
\begin{equation}\label{eq:A.R}
    \hat R = R\hat R R + R^2 -R,
\end{equation}
which follows from the relation $\hat RR+ \hat R+R=0$. Namely, using
\eqref{eq:A.R}, we first show that $\hat R$ is a Hilbert--Schmidt operator,
then that $\hat R\in \fG_2$, and finally that $\hat R\in \fG^1_2$.

Since $R\in \fG^1_2\subset\fS_2$ and $\fS_2$ is an ideal in the algebra of
all bounded operators, we conclude from~\eqref{eq:A.R} that $\hat R
\in\fS_2$.

Next we show that $\hat R$ belongs to $\fG_2$. To this end it suffices to
prove the inclusions $\fG_2\cdot \fS_2 \cdot \fG_2 \subset \fG_2$ and
$\fG_2\cdot \fG_2 \subset \fG_2$. Assume that $R_1,R_3\in \fG_2$, $R_2\in
\fS_2$, and let $R_{ij} := R_iR_j$ and $R_{123}=R_1R_2R_3$. Then the kernel
$r_{12}$ of the operator $R_{12}$ is given by
\[
    r_{12}(x,t) = \int_0^1 r_1(x,s)r_2(s,t)\,ds
        = \bigl(\, r_1(x,\cdot),\overline{r_2(\cdot,t)}\,\bigr)_{L_2(0,1)},
\]
and, using the Cauchy--Schwarz--Bunyakowski inequality, we find that
$r_{12}(x,\cdot)$ belongs to $L_2(0,1)$ and
\[
    \|r_{12}(x,\cdot)\|_{L_2(0,1)}
        \le \|r_1(x,\cdot)\|_{L_2(0,1)} \,
            \|r_2\|_{L_2((0,1)^2)}.
\]
By linearity we also get
\[
    \|r_{12}(x,\cdot) - r_{12}(x',\cdot)\|_{L_2(0,1)}
        \le \|r_1(x,\cdot)- r_1(x',\cdot)\|_{L_2(0,1)} \,
            \|r_2\|_{L_2((0,1)^2)},
\]
which yields continuity of the mapping
 $[0,1]\ni x\mapsto r_{12}(x,\cdot)\in L_2(0,1)$. Similar arguments show
that the kernels $r_{13}$ and $r_{123}$ belong to $G_2$ (in fact, they are
even continuous on $(0,1)^2$). Relation~\ref{eq:A.R} now implies that $\hat
R\in\fG_2$.

Assume now that $R_1,R_3\in\fG^1_2$ and $R_2\in\fG_2$; then
\[
    \frac{dr_{12}(x,t)}{dx} =
        \frac{d}{dx}\bigl(r_1(x,\cdot),\overline{r_2(\cdot,t)}\,\bigr)_{L_2(0,1)}
        =\bigl(r'_1(x,\cdot),\overline{r_2(\cdot,t)}\,\bigr)_{L_2(0,1)}
\]
and thus $r_{12}$ is continuous in~$t$ and once continuously differentiable
in~$x$. Similar arguments show that $R_1R_2R_3\in \fG_2^1$ and
$R_1R_3\in\fG_2^1$. Using this observation in~\eqref{eq:A.R}, we conclude
that $\hat R \in \fG^1_2$, and the lemma is proved.
\end{proof}

\begin{lemma}\label{lem:A.l}
Assume that $\phi\in\Phi\cap W^1_2(0,2)$ and
\begin{equation}\label{eq:A.fact}
    I+F_\phi=(I+K)^{-1}(I+K^\top)^{-1}, \quad K\in \fS_2.
\end{equation}
Then the kernels $k$ and $l$ of the operators $K$ and
$L:=(I+K)^{-1}-I$ respectively have the following property: for
$0\le t\le x\le 1$,
\begin{equation}\label{eq:A.repr}
\begin{aligned}
    k(x,t) & = - f_\phi (x,t) + k_1(x,t),\\
    l(x,t) & =   f_\phi (x,t) + l_1(x,t),\\
\end{aligned}
\end{equation}
where the kernels $k_1$ and $l_1$ are continuously differentiable.
\end{lemma}

\begin{proof}
Applying $I+K$ to both sides of equation~\eqref{eq:A.fact},
rewriting the resulting equality in terms of kernels, and
recalling that $K^\top$ has an upper-diagonal kernel, we arrive at
the so-called Gelfand--Levitan--Marchenko (GLM) equation
\[
    k(x,t) + f_{\phi}(x,t) + \int_0^x k(x,s) f_\phi(s,t)\,ds = 0,
        \qquad 0 \le t \le x \le 1.
\]
Fixing $x\in(0,1)$ in the GLM equation and denoting
\[
    g_x(t):= \begin{cases}
                f_\phi(x,t)& \quad \text{if} \quad 0 < t  <  x < 1,\\
                         0 & \quad \text{if} \quad 0 < x \le t < 1,
                  \end{cases}
\]
we conclude that
\[
    (I+F_\phi) k(x,\cdot) = - g_x.
\]

Since the operator $I+F_\phi$ is invertible and $F_\phi\in\fG^1_2$
by virtue of Lemma~\ref{lem:A.Fphi}, Lemma~\ref{lem:A.invert}
implies that the operator $R:=(I+F_\phi)^{-1}-I$ belongs to
$\fG^1_2$. In particular, with $r$ being the kernel of $R$, we
have for $x>t$
\[
    k(x,t)  = - g_x(t) - \int_0^1 r(t,s)g_x(s)\,ds
            = - f_\phi(x,t) - \int_0^x r(t,s)f_\phi(x,s)\,ds.
\]
Since the function
\[
    k_1(x,t) := - \int_0^x r(t,s)f_\phi(x,s)\,ds
\]
is easily seen to be continuously differentiable on $(0,1)^2$, the
required representation follows.

The operator $L$ satisfies the relation
\[
   L = F_\phi + K^\top + F_\phi K^\top,
\]
or, in terms of kernels,
\[
   l(x,t) = f_\phi(x,t) + \int_0^t f_\phi(x,s) k(t,s)\, ds, \quad
   0\le t\le x\le1.
\]
Now the derived representation~\eqref{eq:A.repr} for the kernel
$k$ implies that the function
\[
   l_1(x,t) := \int_0^t f_\phi(x,s) k(t,s)\, ds
\]
is continuously differentiable on $(0,1)^2$, and the proof is
complete.
\end{proof}

\begin{proof}[Proof of Theorem~\ref{thm:3.smooth}]
Assume that $\si\in W^1_2(0,1)$ and show that then the function $\phi$
of~\eqref{eq:3.phi} belongs to $W^1_2(0,2)$. To this end it suffices to
prove that the kernel $l_\si$ has suitable smoothness properties.

We recall that the function
$a(u,v):=l_\si(\tfrac{u+v}2,\tfrac{u-v}2)$ is continuous on the
set $\ti\Omega^+:=\{(u,v)\mid 0\le v\le u \le 2\}$ and satisfies
there the integral equation
\begin{align*}
        a(u,v) &= -\frac14\int_{v}^{u}d\al\int_0^{v}
                    q\Bigl(\frac{\al-\be}2\Bigr) a(\al,\be)\,d\be\\
                & \quad -\frac12\int_{0}^{v}d\al\int_0^{\al}
                    q\Bigl(\frac{\al-\be}2\Bigr)a(\al,\be)\,d\be\\
                &\quad -\frac12\biggl[\int_0^{u/2}q(\al)\,d\al
                    + \int_0^{v/2}q(\be)\,d\be\biggr]
\end{align*}
with $q:=\si'$. This integral equation implies that the function
\[
        \ti a(u,v):= a(u,v) + \frac12\biggl[\int_0^{u/2}q(\al)\,d\al
                            + \int_0^{v/2}q(\be)\,d\be\biggr]
\]
is continuously differentiable in $\ti\Omega^+$. As a result,
$l_\si$ has the representation
\[
        l_\si(x,t) = \ti a(x+t,x-t)
                     - \tfrac12\si\bigl(\tfrac{x+t}2\bigr)
                     - \tfrac12\si\bigl(\tfrac{x-t}2\bigr) +
                     \si(0),
\]
which implies that the function $\int_0^x l_\si^2(x,t)\,dt$ is in
$W^1_2(0,1)$, so that $\phi\in W^1_2(0,1)$ as well.

Conversely, let $\phi\in \Phi\cap W^1_2(0,2)$ and determine $\si$
through relation~\eqref{eq:3.phi}. In other words, with an
operator $K\in \fS_2^+$ satisfying~\eqref{eq:A.fact} and $l$ being
the kernel of the operator $L:=(I+K)^{-1}-I \in \fS_2^+$, we have
\[
        \si(x) = -2\phi(2x) + 2\int_0^x l^2(x,t)\,dt,
        \quad x\in(0,1).
\]
Using representation~\eqref{eq:A.repr}, one can easily show that
the function
 \(
                \int_0^x l^2(x,t)\,dt
 \)
belongs to $W^1_2[0,1]$. Henceforth $\si\in W^1_2(0,1)$, and the
proof is complete.
\end{proof}


\end{document}